\documentclass{amsart}

\usepackage[colorlinks,urlcolor=blue,citecolor=blue,filecolor=blue,linkcolor=blue]{hyperref}

\usepackage[applemac]{inputenc}
\usepackage[english]{babel}
\usepackage{amssymb}
\usepackage{amsfonts}
\usepackage{amscd}
\usepackage{graphicx}
\usepackage{epstopdf}
\usepackage{comment}

\newcommand{\translation}[1]{$[\![$#1$]\!]$}

\newcommand{\R}{{\mathbb{R}}}
\newcommand{\G}{{\mathbb{G}}}

\newcommand{\W}{{\mathbb{W}}}
\newcommand{\V}{{\mathbb{V}}}
\newcommand{\Pu}{{\mathbb{P}}}

\newcommand{\Lc}{{\mathcal{L}}}
\newcommand{\Vc}{{\mathcal{V}}} 
\newcommand{\Pc}{{\mathcal{P}}} 

\newtheorem{theorem}{Theorem}[section]

\newtheorem{proposition}[theorem]{Proposition}

\newtheorem{metatheorem}[theorem]{Metatheorem}

\theoremstyle{definition}
\theoremstyle{remark}
\numberwithin{equation}{section}

\newcommand{\achtung}{\marginpar{**********}}
\newcommand{\pq}[1]{[\![#1]\!]}

\newcommand{\labelpag}[1]{\refstepcounter{enumi}\label{#1}}

\begin{document}
\title[Peano on  surface  area]{Peano on definition of surface  area}
\author{Gabriele H. Greco}
\address{Independent Researcher, Professor of Mathematical Analysis,
Trento, Italy}
\email{gabriele.h.greco@gmail.com, https://sites.google.com/site/ghgreco/}

\author{Sonia Mazzucchi}
\address{Dipartimento di Matematica\\
Universit\`{a} di Trento, 38123 Povo (TN), Italy}
\email{sonia.mazzucchi@unitn.it}

\author{Enrico M. Pagani}
\address{Dipartimento di Matematica\\
Universit\`{a} di Trento, 38123 Povo (TN), Italy}
\email{enrico.pagani@unitn.it}

\dedicatory{On the occasion of the $150^{th}$ anniversary of the birth of Giuseppe Peano}
\date{%
December 8, 2014
}

\begin{abstract}
In this paper we investigate the evolution of the concept of area
in \textsc{Peano}'s works, taking into account the main role played 
by \textsc{Grassmann}'s geometric-vector calculus
and \textsc{Peano}'s theory on derivative of measures.
\emph{Geometric} (1887) and \emph{bi-vectorial} (1888) \textsc{Pea\-no}'s approaches 
to surface area mark the development of 
this topic during the first half of the last century. 
In the sequel we will present some significative contributions on surface area that 
are inspired and/or closely related to \textsc{Peano}'s definition. 
\end{abstract}

\maketitle

\section{Introduction}

In 1882 \textsc{Peano} at the age of 24 discovers that the definition of area of a surface presented by \textsc{Serret} in his \emph{Course d'\,Analyse} \cite[(1868) vol.\,2, p.\,296]{serret} was not correct. According to \textsc{Serret}'s proposal, the area of a surface should be given by the limit 
of the area of the inscribed polyhedral surfaces, but  this definition cannot be applied even to a cylindrical surface. In fact \textsc{Peano} observes that in this case it is possible to choose a suitable sequence of inscribed polyhedral surfaces whose areas converge to infinity (see \cite[(1890)]{peano_area1890}, \cite[(1902) p. 300-301]{formulaire4}).

\textsc{Genocchi}, \textsc{Peano}'s teacher, dampens the enthusiasm of the young mathematician, by communicating him that a similar counterexample was already been discovered by \textsc{Schwarz}. 
In fact, in a letter of May 26, 1882, \textsc{Genocchi} writes \textsc{Schwarz} \cite[(1890) vol.\,2, p.\,369]{schwarz}:
\begin{quotation}
 C'est précisément Mr.\,Peano, qui m'amène à vous parler d'un autre sujet. Devant aborder la quadrature des surfaces courbes, il s'est aper\c cu que  la définition d'une aire courbe donnée par Serret n'était pas bonne, et m'a expliqué les raisons qui ne lui permettaient pas de l'adopter. Alors je l'ai informé du jugement que vous en aviez porté dans plusieurs de vos lettres (20 et 26 décembre 1880, 8 janvier 1881), ce qui l'a beaucoup intéressé.
\end{quotation}
\textsc{Genocchi} and \textsc{Schwarz}\,\footnote{ \ \textsc{Schwarz} communicates his counterexample also to \textsc{Casorati} and to \textsc{Beltrami} (1880); 
see the correspondence between \textsc{Casorati} and \textsc{Peano} in \textsc{Gabba} \cite[(1957)]{gabba_carteggio1957}.}
were conscious of the problem and of the lack of a ``correct'' definition of surface area, suitable to handle at least the area of elementary figures. 
In 1882 \textsc{Genocchi} writes another letter to \textsc{Schwarz} and invites him to propose an alternative definition, but \textsc{Schwarz} declines and stresses the difficulties:
\begin{quote}
Vous avez voulu que je donne la rectification de la d\'efinition incompl\`ete; mais ce n'est pas facile. On peut rectifier cette d\'efinition {\it de plusieures mani\`eres} et il me semble qu'il suffit de donner express\'ement une seule possibilit\'e qui convient avec la d\'efinition donn\'ee par Sturm.
\end{quote}
\textsc{Sturm}'s definition \cite[(1877) vol.\,1, p.\,427]{sturm} is the following:
\begin{quotation}
On appelle aire d'une surface courbe, terminée à un contour quelconque, la limite vers laquelle tend l'aire d'une surface polyédrique composée de faces planes, qui en diminuant toutes indéfiniment, tendent à devenir tangentes à la surface considérée. On suppose d'ailleurs que le contour qui termine la surface polyédrique se rapproche indéfiniment de celui qui termine la surface courbe.
\end{quotation}

We can think that the drawbacks communicated by \textsc{Schwarz} to \textsc{Genocchi} were also 
due to the lack of a choice criterion between the several possible definitions of surface area. 
In any case a good definition of area should at least be compatible with the Lagrange
formula of area of a Cartesian surface \footnote{ \ Nowadays we know that the Lagrange formula is 
sufficient to define the area of a  $C^1$-submanifold, but the extension of the formula (\ref{Lagrange})
from a rectangle $D\/$ to a more general $2\/$-dimensional set is not trivial and  hides some pitfalls.}: 
\begin{equation} \label{Lagrange}
\iint_D  \sqrt{1+|\nabla f(x,y)|^2} \, dx dy
\end{equation} 
for $C^1$ functions $f:D\to\R$ ($D$ being any rectangular subset of $\R^2$).

\textsc{Peano}'s definition of area in \emph{Applicazioni geometriche del calcolo infinitesimale} \cite[(1887)\,p.\,164]{peano87} overcomes the drawbacks of \textsc{Serret}'s approach, 
yielding the Lagrange formula \eqref{Lagrange}.
\textsc{Peano}'s proposal is deeply influenced by \textsc{Grassmann}'s geometric-vector calculus in affine spaces, that gives a mathematical formalization of geometrical and physical concepts (vectors, pair of vectors, moment and so on) and allows also to take into account properties related to orientation, without using drawings or tricky and intuitive constructions. 
It is not surprising that \textsc{Peano}'s definition via \textsc{Grassmann}'s calculus is suitable 
to handle oriented integrals and, consequently, to prove main results (such as Stokes theorem and Green formula), and to develop formulae leading to the integration of 2-forms of \textsc{Cartan} 
\footnote{ \ In  1899 \textsc{Cartan} \cite{cartan1899} introduced the calculus of differential forms 
\begin{quote}
Ce calcul pr\'esente aussi de nombreuses analogie avec le calcul de Grassmann; il est d'ailleures identique au calcul g\'eom\'etrique dont se sert M. Burali-Forti dans un Livre r\'ecent (\emph{Introduction 
\`a la G\'eom\'etrie diff\'erentielle, suivant la m\'ethode de Grassmann}, Gauthier-Villars, 1898).
\end{quote}
\textsc{Burali-Forti} was one of the prominent scholars of \textsc{Peano}. Together with
\textsc{Marcolongo} he developed and applied \textsc{Grassmann}'s vector calculus to
geometry, mechanics and physics.   
}.

Besides \textsc{Peano}'s proposal, in the literature several definitions of surface area have been given 
\footnote{ \ For a detailed presentation of the several possible definitions of surface area see 
\textsc{Cesari} \cite[(1954)]{cesari-book} and
\textsc{Federer} \cite[(1969)]{federer}}: nowadays the most famous and commonly accepted as 
definitive are grounded on Hausdorff measures.  

The aim of the present paper is the investigation of the evolution and use of the concept of area
in \textsc{Peano}'s works, 
taking into account the main role played by \textsc{Grassmann}'s geometric-vector calculus
and \textsc{Peano}'s theory on derivative of measures.
\textsc{Peano}'s approach to surface area marks the development of this topic during the
first half of the last century. 
In the sequel we will present contributions concerning surface area that are inspired and/or 
closely related to \textsc{Peano}'s definition. 


\textsc{Peano}'s definition of measure of surfaces is grounded 
on elementary formulae of area of planar polygons  
(see Theorems \ref{2.1}, \ref{2.2} and their proofs).
The surprising absence of results concerning area of planar polygons in several modern encyclopedic books 
(see for example \textsc{Alexandrov} \cite[(2005)]{alexandrov} and \textsc{Berger} \cite[(1977)]{berger}) 
motivates us to try to trace the history of such formulae that, as we shall see, are deeply connected with 
statics and can be found in their final form in the works by \textsc{M\"obius} and \textsc{Bellavitis}.  
The generalization of the formula of area from planar to non-planar polygons and to
closed curves allowed \textsc{Peano} to specify and to evaluate area of surface at an infinitesimal level
(see Section \ref{sez-Peano-area}).

The paper is organized as follows. In Section \ref{sez-vec-calculus} the main definitions and results on Grassmann's geometric-vector calculus are presented in a modern fashion, according to 
\textsc{Greco, Pagani} \cite[(2010)]{Greco-Pagani}. 
In Section \ref{sez-Mobius-Bell} the historical development of the formulae of area of polygons and 
volume of polyhedra is investigated. 
Section \ref{sez-Peano-area} is devoted to the description of \textsc{Peano}'s definition of area.
In Section \ref{bi-vectors} we recall \textsc{Peano}'s bi-vector integral formula and other
ways of associating a \emph{number} to a given oriented closed curve.
In Section \ref{uso} we will list main propositions and theorems about area 
given by \textsc{Peano} in his works.
In Section \ref{Peano-influenza} we recall some significant mathematical contributions
inspired and/or closely related to \textsc{Peano}'s definition. In particular 
we present the re-formulations of \textsc{Peano}'s and \textsc{Ge\"ocze}'s area 
(due to several mathematicians) 
in order to make them coincident with \textsc{Lebesgue}'s area.
%
%
%
%
%
%

This article concerns some historical aspects. From a methodological point of view, we are focussed on primary sources, that is on mathematical facts and not on the elaborations or interpretations of these fact by other Scholars of history of mathematics.


\section{Grassmann-Peano geometric-vector calculus \\ on three dimensional affine spaces}\label{sez-vec-calculus}

We present here the Grassmann-Peano geometric-vector calculus, as described in \textsc{Greco, Pagani}
\cite[(2010)]{Greco-Pagani}. The aim of this section is to understand the mathematical basis used by \textsc{Peano} in the construction of his notion of area. Such a formalism will be useful not only to clarify the genesis of  the vectorial formulae for area of polygons and volume of polyhedra, but also
to understand the deep connection between some concepts of statics (points, applied forces,
momenta, Poinsot pairs and so on) and of geometry (geometric forms of first, second and $n\/$-degree,
namely, points, vectors, bi-points, tri-points, quadri-points, and so on).

\textsc{Peano} is one of the first mathematician who presents \textsc{Grassmann}'s work 
\cite[(1844)]{grassmann1844}, \cite[(1862)]{grassmann1862}, \cite{grassmannW} to the mathematical community. 
Actually he rebuilds \textsc{Grassmann}'s calculus using a original ``functional'' approach
that relies only on the assignment of a volume form on a given affine space (see  \textsc{Greco, Pagani}
\cite{Greco-Pagani} for a detailed presentation of this subject).

For convenience of the reader we choose to rebuild here \emph{Grassmann graded exterior algebra 
on an affine space} (Grassmann affine algebra, for short),
using an approach based on the usual notion of graded exterior algebra on a vector space.

The starting point for the construction of Grassmann affine algebra 
on the ordinary $3$-dimensional affine space is the introduction of a $4$-dimensional \emph{M\"obius space}, i.e., a couple $(\W,\omega)$, where $\W$ is $4$-dimensional vector space  
and $\omega:\W\to\R$ is  a non-vanishing linear form,  called \emph{mass}.
Given the M\"obius space $(\W,\omega)$, let us consider the subspace 
$\V:=\{w \in \W : \omega(w)=0\}$ of $\W$ and the subset $\Pu:=\{w \in \W : \omega(w)=1\}$. 
Elements of $\V$ and $\Pu$ will be called $\omega$-{\it vectors} and $\omega$-{\it points}
of $(\W,\omega)$, respectively.

The affine space $(\Pu, \V, -)$ (where $-$ stands for the difference between elements of $\W$) 
may be identified with the 3-dimensional Euclidean space. Therefore, the form $\omega\/$ allows a non ambiguous selection of the $\omega$-\emph{vectors} and the $\omega$-\emph{points} from the elements of the M\"obius space 
$\W\/$.
The elements of $\W\/$ with $\omega (x)\neq 0$ are called \emph{weighted $\omega$-points}.	

Let us consider the graded exterior algebras $\G(\W)$ and $\G(\V)$ on the
vector spaces $\W\/$ and $\V\/$ respectively. We have explicitly 
$$\G(\W)=\Lambda^0(\W)\oplus\Lambda^1(\W)\oplus\Lambda^2(\W)\oplus\Lambda^3(\W)\oplus\Lambda^4(\W)$$
where $\Lambda^0(\W):=\R$ and $\Lambda^k(\W)$, $k=1,\dots,  4$, is the vector space generated by the
products of $k$ vectors of $\W\/$. The elements of $\Lambda^k(\W)$ are called \emph{geometric forms of degree} $k$. 
Since $\W\/$ is the vector space generated by
$\omega$-points, it is worth observing that $\Lambda^k(\W)$
is generated by the products of $k\/$ $\omega$-points.
 In a similar way we have
$$\G(\V)=\Lambda^0(\V)\oplus\Lambda^1(\V)\oplus\Lambda^2(\V)\oplus\Lambda^3(\V)$$ 
where $\Lambda^0(\V):=\R$ and $\Lambda^k(\V)$, $k=1, 2, 3$, are linear combinations of products of $k$ $\omega$-vectors. 
The elements of $\Lambda^k(\V)$ are called \emph{geometric vector forms of degree} $k$.
\footnote{ \ 
The product in $\G(\W)$ and in $\G(\V)$ will be denoted  as juxtaposition of symbols. 
Recall that the algebras $\G(\W)$ and $\G(\V)$ are anticommutative, i.e. $xy= (-1)^{rs}yx$, for any $x\in\Lambda^r(\W),y\in\Lambda^s(\W)$. Clearly, $\G(\V)$ is a sub-algebra  of $\G(\W)$;
moreover, due to anti-commutativity of the product and to the 
dimension of the space $\W$ and $\V\/$ we have that $\Lambda^k(\W) = {0}$ (for $k > 4$) and
$\Lambda^k(\V) = {0}$ (for $k > 3$).}


The linear form $\omega:\W \to \R$ can be extended to a linear map from the whole $\G(\W)$ to $\G(\V)$ by means of the following relations
$$ \omega (1)=0$$
$$\omega(P_0)=1$$
$$\omega(P_0P_1)=P_1-P_0$$
$$\omega(P_0P_1P_2)=(P_1-P_0)(P_2-P_0)$$
$$\omega(P_0P_1P_2P_3)=(P_1-P_0)(P_2-P_0)(P_3-P_0),$$
for every $P_0, P_1, P_2, P_3 \in \Pu.$
\footnote{ \ 
The extension of the linear form $\omega$ on the whole graded algebra $\G(\W)$ can be uniquely determined by the conditions
\begin{enumerate}
\item $\omega:\Lambda^k(\W)\to\Lambda^{k-1}(\V)$ (we assume $\Lambda^{-1}(\V):=\{0\}$);
\item $\omega(xy)=\omega (x)y+(-1)^{{\rm deg}(x)}x\omega (y)$ (\emph{graded Leibnitz rule}).
\end{enumerate}
}
Actually  $\Lambda^k(\V)$ is a vector subspace of $\Lambda^k(\W)$ and the linear map $\omega $ connects the graded algebras $\G(\W)$ and $\G(\V)$ in the following way:
\begin{equation} \label{triviality}
\omega (\Lambda^k(\W))=\Lambda^{k-1}(\V)={\rm Ker }\,(\omega_{|\Lambda^{k-1}(\W)}) \, \, k = 0, \cdots, 4. 
\end{equation}  
Restrictions of $\omega$ to $\Lambda^k(\W)$, denoted by $\omega_k$, are called $(k-1)$-\emph{vector-masses} because,
by the first equality of formula \eqref{triviality}, $\omega$ transforms a $k$-degree geometric form
into geometric vector forms of  $(k-1)$-degree. The second equality says that a $k$-vector-mass is 
null on geometric vector forms of degree $k$;
in particular $\omega \circ \omega=0$ and the \emph{reduction formula} holds:
\begin{equation}
\label{riduzione}
x=P\omega (x)+\omega (Px),\qquad \forall P\in \Pu,\;x\in \Lambda ^k(\W) \,, k=0,1,2,3,4.
\end{equation}

A quadri-point $ABCD$ (with $A,B,C,D \in \Pu\/$ regarded as the four vertices of a tetrahedron) suggests the construction of particular basis of $\Lambda^r(\W)$, whenever they are not co-planar  (i.e. the vectors $B-A,C-A,D-A$ are linearly independent). The four vertices $A,B,C,D$ are a basis of $\Lambda^1(\W)$, the six bi-points $AB,AC,AD, BC,BD,CD$, corresponding to the six edges of the tetrahedron,  are a basis of $\Lambda^2(\W)$, the four tri-points $ABC,ACD,ABD, BCD$, corresponding to the four faces of the tetrahedron,  are a basis of 
$\Lambda^3(\W)$, and  the quadri-point $ABCD$ is a basis of $\Lambda^4(\W)$.

If $A',B',C',D'$ are the four vertices of another tetrahedron,  then 
\begin{equation}\label{det}A'B'C'D'=\frac{\det (B'-A',C'-A',D'-A')}{\det (B-A,C-A,D-A)}ABCD,
\end{equation}
 where $(B'-A',C'-A',D'-A')$ is the $3\times 3 $  matrix whose columns are the coordinates of the vectors $B'-A',C'-A',D'-A'$ along the basis $B-A,C-A,D-A$. Equality \eqref{det} enlightens the geometrical interpretation of a quadri-point in terms of an oriented volume.  
The equality between two elements $x,y \in \Lambda^r(\W)$ can be expressed by means the following condition:
\begin{equation}\label{equality}x=y  \qquad \Leftrightarrow \qquad   xz=yz \quad \forall z\in \Lambda^{4-r}(\W) 
\end{equation}

Several elements of the graded exterior algebras $\G (\W) $ and  $\G (\V) $ admit interesting geometrical and mechanical interpretation. 

Let $A,B,C \in \Pu\/$. The bi-point $AB\/$ can be seen as the applied vector $A(B-A)$ 
(for instance, as a ``force'' $B-A\/$ ``applied'' in $A$), and the tri-point 
$ABC\/$ can be represented by a triangle or by an applied bi-vector $A(B-A)(C-A)$ (applied in $A$). 
Elements of $ \Lambda^2 (\V) $, i.e. the bi-vectors, can be seen as Poinsot couples or as  oriented boundary of triangles, and  the elements of $ \Lambda^3 (\V) $, i.e.  tri-vectors, as oriented surfaces of tetrahedrons.

Besides mechanical interpretations, a system of applied forces can be represented by an element of 
$\Lambda^2(\W)$, more precisely as a sum of bi-points. 
The equivalence between  two systems of
applied forces $\{A_iB_i\}_{i=1, \dots, n}$ and $\{C_jD_j\}_{j= 1, \dots , m}$
can be expressed as the equality between the corresponding elements 
$\sum_i A_iB_i$ and $\sum_j C_jD_j$ of  $\Lambda^2(\W)$. 
Indeed by means of equation \eqref{det}, given a bi-point $PQ$, the product $A_iB_iPQ$  can be recognized as the axial moment of the force $A_iB_i$ with respect to the axis passing through $P$ and $Q$. 
As a consequence the equality \eqref{equality} between elements of $\Lambda^2(\W)$ reduces
the equivalence between two systems of applied forces to the equality  
of their axial moments with respect to every axis.
As a particular case, a system of forces with vanishing resultant (a Poinsot couple), can be represented by an element of 
$\Lambda^2(\V)$. 
It is interesting to note that \textsc{Poinsot}'s theorem concerning the sum of Poinsot couples
emerges naturally from the structure of vector space of $\Lambda^2(\V)$.
 
Pursuing the analogy with statics, the image of the  operator $\omega$ acting on $\Lambda^2(\W)$ represents the resultant of a system of forces (a special case of the $1$-vector-mass introduced above). The reduction formula \eqref{riduzione} can be directly translated into the reduction formula for a system of forces: given an arbitrary point $P$, a system of forces $x$ is equivalent to a system formed by the resultant $\omega (x)$ applied in $P$, and by the Poinsot couple $\omega (Px)$.





The formalism presented so far, allows a direct proof of the following results:
\begin{theorem}[Area of a plane polygon%
\,\footnote{\ See \textsc{Peano}'s 
\emph{Applicazioni geometriche del calcolo infinitesimale} {\cite[(1887)\,p.\,21]{peano87}},
\emph{Calcolo geometrico} {\cite[(1888)\,p.\,59]{peano88}},
\emph{Lezioni di analisi infinitesimale} {\cite[(1893)\,vol.\,II,\,p.\,32]{peano1893}}.}] \label{2.1}
For any planar polygon with consecutive vertices $A_1,\dots , A_{n}$, the sum 
\begin{equation}\label{area-poligoni}
\sum_{i=1}^{n} PA_iA_{i+1}
\end{equation}
(with $A_{n+1} = A_1$) does not depend on the choice of the point $P$, with $P$ belonging to the plane
of the vertices.
\end{theorem}

The vector space of third degree forms in a plane is $1$-dimensional; a base for this space is 
provided by an arbitrary triangle $RST$ with non collinear vertices. 
Then $PA_iA_{i+1} = a \, RST$, where $a$ is the oriented area of $PA_iA_{i+1}$ with
respect to $RST$ \footnote{ \ 
In other words, $|a|$ is the ratio between the areas of the two triangles; the sign of $a$
is positive if the triangles have the same orientation.}. Therefore formula (\ref{area-poligoni})
gives the sum of the oriented areas of triangles $PA_iA_{i+1}$, termed by \textsc{Peano} as
the \emph{area bounded by the oriented closed polygonal line} $A_1,\dots , A_{n}, A_{n+1}= A_1$.
As observed by \textsc{Peano}, this area coincides with the usual measure of area if
the polygonal line is convex or, more generally, is not interlaced. 

\begin{theorem}[Area of a non-planar polygon%
\,\footnote{\ See \textsc{Peano}'s 
\emph{Calcolo geometrico} {\cite[(1888)\,p.\,137]{peano88}}.}] \label{2.2}
For any closed polygonal line (not necessarily planar) there exists a triangle such that the area of any projection 
of the polygonal line on an arbitrary plane is equal to the area of the 
projection of the triangle.~\footnote{\ In this case the polygon is said ``equipollent'' to the triangle.}
\end{theorem}
\begin{theorem}[Volume of an oriented polyhedron%
\,\footnote{\ See \textsc{Peano}'s
\emph{Applicazioni geometriche del calcolo infinitesimale} {\cite[(1887)\,p.\,26-27]{peano87}},
\emph{Calcolo geometrico} {\cite[(1888)\,p.\,66]{peano88}},
\emph{Lezioni di analisi infinitesimale} {\cite[(1893)\,vol.\,II,\,p.\,35]{peano1893}}.}]  \label{2.3}
Let us consider a closed oriented polyhedral surface made of triangular faces $A_iB_iC_i$, $i=1\dots , n$. The sum of the oriented volumes 
\begin{equation}\label{volume-poliedri}
\sum_{i=1}^{n} PA_iB_iC_i
\end{equation}
of the tetrahedra $PA_iB_iC_i$  
does not depend on the choice of the vertex $P$.
\end{theorem}

For convenience of the reader, we give the proofs of the previous theorems accordingly to Section \ref{sez-vec-calculus} 
on Grassmann-Peano geometric-vector calculus.

\begin{proof}[{\bf Proof of theorem \ref{2.1}}]
Denote by $\pi$ the plane of the vertices $A_1, \dots, A_n$.
Let us consider the element $x\in\Lambda^2(\W)$ given by $x=\sum_{i=1}^nA_i{A_{i+1}}$. 
Since $\omega_2 (x) = \sum_{i=1}^n(A_{i+1} - A_i)= A_{n+1} - A_1 =0$, we have
$x \in \ker (\omega_2)=\Lambda^2(\V)$; therefore, $x$ is a bi-vector and, hence, there exist three points $X,Y,Z$ in 
the plane $\pi$ such that $x=(Y-X)(Z-X)$. Given a generic point $P$ in the plane $\pi$, 
we have $\omega_4(PXYZ)= (X-P)(Y-P)(Z-P)=0$ (the three vectors are linearly dependent); therefore,
by the reduction formula \eqref{riduzione}, we have 
$$XYZ=P\omega_3 (XYZ)+\omega_4 (PXYZ)= P\omega_3 (XYZ)=Px= \sum_{i=1}^n P A_i{A_{i+1}}$$
and the conclusion follows.
\end{proof}

\begin{proof}[{\bf Proof of theorem \ref{2.2}}]
Let us consider a (non planar) polygons with vertex $A_1, $ $\dots , $ $A_n$, with $A_{n+1}=A_1$. By reasoning as in the proof 
of Theorem \ref{2.1}, the element $x\in\Lambda^2(\W)$ given by $x=\sum_{i=1}^nA_i{A_{i+1}}$ can be represented by a 
bi-vector and there exist three points $X,Y,Z$ such that $x=(Y-X)(Z-X)$.
Given a generic plane $\gamma$ and a generic point $A$, let us consider a unitary vector $u$ orthogonal to the 
plane; consequently $Au$ represents a ``vector applied'' applied at the point $A$. 
Relation \eqref{equality} implies that 
$(Y-X)(Z-X)Au=\sum_{i=1}^nA_i{A_{i+1}}Au$. By formula (\ref{det}) the 
element $(Y-X)(Z-X)Au\in \Lambda ^4(\W)$ represents the volume of the tetrahedron with unitary height and area of 
basis equal to the area of the projection of the triangle $XYZ$ on the plane $\gamma$, and the conclusion follows.
\end{proof}

\begin{proof}[{\bf Proof of theorem \ref{2.3}}]
The closedness of the polyhedral surface made of triangular faces $A_iB_iC_i$  amount to the condition 
$\omega_3(\sum_{i=1}^n A_iB_iC_i)=0$, which implies that the element $\sum_{i=1}^n A_iB_iC_i\in \Lambda ^3(\W)$ 
is a tri-vector. 
This implies that there exist four points $X,Y,Z,U$, such that $\sum_{i=1}^n A_iB_iC_i=(Y-X)(Z-X)(U-X)$. 
For any point $P$, the reduction formula \eqref{riduzione} gives 
\begin{eqnarray} \label{vol3}
XYZU &=&P\omega_4(XYZU) +\omega_5(PXYZU) = P\omega_4(XYZU)  \cr 
&=& P(Y-X)(Z-X)(U-X)=\sum_{i=1}^n PA_iB_iC_i \, \, \,,
\end{eqnarray}
and the conclusion follows.
\footnote{\ The independence of the sum of volumes (\ref{vol3}) 
has been proved starting from the equality $\omega_3(\sum_{i=1}^n A_iB_iC_i)=0$.
It is worth noting that the converse is still valid; in other words the equality $\omega_3(\sum_{i=1}^n A_iB_iC_i)=0$ holds if and only if  ``$v\sum_{i=1}^n A_iB_iC_i=0$ for every $v\in\V$''.
}
\end{proof}

%

\section{M\"obius and Bellavitis on the area of polygons \\ and volume of polyhedra}\label{sez-Mobius-Bell}

At the beginning of the $19^{th}$ century an increasing interest is devoted to the study of polygons and polyhedra. This interest is paved by the researches by \textsc{Legendre} and \textsc{Poinsot}, who follow the way traced by \textsc{Euclid}, \textsc{Kepler}, \textsc{Descartes} and \textsc{Euler}. \textsc{Legendre} in 1794 gives a proof of the famous \textsc{Euler}'s formula (1750) for polyhedra:
\begin{equation}\label{eulero}V-E+F=2,\end{equation}
where $V, E$ and $F$ denote the number of vertices, edges and faces, respectively.
On the other hand, \textsc{Poinsot} \cite{poinsot1810}, according to the ``G\'eom\'etrie de situation'' of \textsc{Leibnitz} \cite{Leibniz},
in 1810 started the classification of polygons and polyhedra, discovering some new ``star polyhedra''. 
In 1813 \textsc{Cauchy} \cite{cauchy1813} gave a new proof of \textsc{Euler}'s formula \eqref{eulero} showing that there are no star polyhedra different from those described by \textsc{Poinsot}. Moreover, urged by \textsc{Legendre}, \textsc{Cauchy} gave the famous rigidity theorem for convex polyhedra, as he said in 
{\em  Sur les polygones et les poly\`edres} \cite[p.\,87]{cauchy1813a}:
\begin{quote}
[...] chercher la d\'emonstration du th\'eor\`eme renferm\'e dans la d\'efinition 9, plac\'ee à la t\`ete du onzi\`eme Livres {\em Elements d'Euclide}, savoir que deux poly\`edres convexes sont \'egaux lorsqu'ils sont compris sous un m\^eme nombre de faces \'egales chacune à chacune.
\end{quote}
%
One of the first book devoted to polyhedra was written by \textsc{Descartes} \cite{Cartesio}, but many other authors devote their efforts to the study of this topic.



An evidence of the importance which was given to polygons and polyhedra in the $19^{th}$ century is the \emph{Gran Prix} ``Perfectionner dans quelque point important la th\'eorie g\'eom\'etrique des poly\`edres'' organized  in 1858 by the Accademy of Sciences of Paris. Indeed  the Accademy decided 
to assign a prize only in presence of a significative and revolutionary contribution to the theory of polyhedra. Several important scientists participate, including \textsc{M\"obius}.
As other participants, \textsc{M\"obius}'s goal was to provide a complete classification of polyhedra, but very 
soon he discovered that this is really an arduous task and decided to change his aims, proposing an innovative work concerning the concept of orientation. 
Despite of this, the Accademy does not judge any contribution sufficiently important and does not assign 
the prize to any participant. 

Among several results present in the mathematical literature, we restrict ourselves to analyze in details the works of \textsc{M\"obius} and \textsc{Bellavitis}, due to their influence on \textsc{Peano}. 
The formula of area of polygons \eqref{area-poligoni} can be found for the first time in 
\textsc{M\"obius}'s {\em Barycentrische Calcul} \cite[(1827)\,p.\,201]{Moebius}, where it appears in a remark, as an application of the analogous formula for triangles and as a direct geometrical consequence of the notion of barycentric coordinates. \textsc{Bellavitis} presents the formula in \emph{Teoremi generali per determinare le aree dei poligoni e i volumi di
poliedri}
\cite[(1834)]{bellavitis1834} as a ``trivial consequence'' of a well known properties due to \textsc{Poinsot}.
\textsc{Bellavitis} says indeed:
%
%
%
%
%
%
\begin{quote}
[La formula dell'area]
esprime la proprietà di un sistema di forze aventi la risultante nulla di produrre una stessa coppia ({\it couple}), in qualunque punto comune tutte esse si trasportino. \footnote{ \pq{%
[In the formula of area] one can see the property satisfied by a system of applied forces with vanishing resultant
to be equivalent to a couple, independently of the common point in which the forces are translated.%
}}
\end{quote} 
Later \textsc{M\"obius} himself deduced the formula \eqref{area-poligoni} of previous Section, as an application 
of statics, in  his book
{\em Der Statik} \cite[(1837)\,p.\,61-64]{moebius1837}.
In our opinion this correlation with statics, where the couples of consecutive vertices ($=$ bi-points) 
of a closed polygonal line are interpreted as forces with vanishing resultant,
is important from a historical point of view and may be emphasized into the following:
\begin{metatheorem} \label{metathm}
The following two propositions are equivalent:
\begin{enumerate}
\item Formula of area \eqref{area-poligoni} for planar polygons holds.
\item Any system of planar forces with vanishing resultant is equivalent to a couple.
\end{enumerate}
\end{metatheorem}
According to \textsc{Bellavitis} also the formula of volume of polyhedra \eqref{volume-poliedri} can be seen as a consequence of the static theorem of \textsc{Poinsot}: ``the sum of couples is a couple''.
Later, references to   formulae \eqref{area-poligoni} and \eqref{volume-poliedri} can be found in
\textsc{Bellavitis}'s {\em Metodo delle equipollenze} \cite[(1838)\,pp.\,95-97]{bellavitis1838b}
and \emph{Sposizione del metodo delle equipollenze} \cite[(1854)]{bellavitis1854}.
\par
Concerning \textsc{M\"obius}, both formulae \eqref{area-poligoni} and \eqref{volume-poliedri} can be found in his article appeared in 1865 {\em \"Uber die Bestimmung des Inhaltes eines Polyhedres} 
\cite[pp.\,486, 494]{moebius1865}. 

The methods of proof of \textsc{Bellavitis} and \textsc{M\"obius} are quite different. \textsc{Bellavitis} is one of the first mathematicians developing vector calculus, and he uses it in most of his proofs. Moreover the deep connection between statics and geometry is strongly emphasized. It is also worthwhile to note that \textsc{Bellavitis} applies the duality relation between polygons and polyhedra, then it is not surprising that both formulae \eqref{area-poligoni} and \eqref{volume-poliedri} appear in the same article.
In the work by \textsc{M\"obius} emerges the revolutionary concept of orientation. 
\textsc{M\"obius} is conscious that orientability of polyhedra is an important condition for the validity of the formula of volume, and it cannot be ignored, as well as he was aware of the existence of non oriented polyhedra. 
In \textsc{Bellavitis}'s work the necessity of orientability is not transparent, and he handles only with polyhedra which are dual of polygons that are oriented by construction. 

We did not find any trace (but we cannot exclude it) concerning area of non-planar polygons either 
in \textsc{M\"obius} or in \textsc{Bellavitis} even if a statement similar to Metatheorem \ref{metathm} 
is still valid for non planar polygons and non planar forces:
\begin{metatheorem} \label{metathm2}
The following two propositions are equivalent:
\begin{enumerate}
\item  For any closed polygonal line (not necessarily planar) there exists a triangle such that the area of any projection 
of the polygonal line on an arbitrary plane is equal to the area of the 
projection of the triangle.
\item Any system of forces with vanishing resultant is equivalent to a couple.
\end{enumerate}
\end{metatheorem}



\section{\textsc{Peano}'s definitions of area}\label{sez-Peano-area}
%
%
%
%
%
%
In \textsc{Peano}'s works we recognize two definitions of area of a non-planar surface,
the first one, referred as \emph{geometric}, is based on the notion of Jordan-Peano area of planar sets;
the second one, referred as \emph{bi-vectorial}, is based on the notion of bi-vector associated to
a closed curve bounding a pieces of surface.
%
%
%
%
%
%

In \emph{Applicazioni Geometriche del Calcolo Infinitesimale} \cite[p.\,164]{peano87} 
\textsc{Peano} introduces his geometric definition of area in the following terms:
\begin{quote}
{\sc Aree di superficie non piane.}
Abbiasi una superficie qualunque. Proiettandola ortogonalmente sopra un piano avremo una figura pia\-na; supporremo che questo abbia un'area propria, e che la superficie data si possa decomporre in parti 
che godano della stessa propriet\`a.\\
Si scomponga la superficie data in parti e, dopo averle trasportate comunque nello spazio, si proiettino queste ortogonalmente su d'uno stesso piano. La somma delle aree di queste proiezioni sarà un'area piana, variabile col variare del modo di divisione della superficie e del modo con cui si dispongono queste parti. Il limite superiore dei valori di quest'area piana si dirà l'{\it area} della superficie data.\\
Si deduce immediatamente dalla definizione che l'area di una superficie qualunque \`e maggiore della sua proiezione ortogonale su d'un piano qualunque.\,%
\footnote{\ 
\pq{\textsc{Area of non planar surfaces}. Let us consider an arbitrary surface. Performing an orthogonal projection on a plane, we get
a plane figure; we assume that this figure have an ``area propria'' (i.e., it is Peano-Jordan measurable)
and that the given surface can be decomposed into parts having the same property.  

Let us decompose the given surface into pieces and, after carrying these pieces arbitrarily in the space, let us project all these pieces on the same plane. The sum of the areas of these projections is a planar area,  
depending on the decomposition of the surface and on the way its pieces are located. 
The supremum of the values of these planar areas will be defined as the {\it area} of the surface.

It follows immediately from this definition that the area of an arbitrary surface is greater 
than its orthogonal projection on an arbitrary plane.}
}
\end{quote}


Paraphrasing the content of the paper \emph{Sulla definizione dell'area di una superficie} 
\cite[(1890), p.\,56] {peano_area1890}, we may have the following bi-vectorial definition of area,
that may help the reader in comparing the two definitions of area given by \textsc{Peano}:

\begin{quote}
Let us consider an arbitrary surface delimited by a closed oriented curve. 
Performing an orthogonal projection on a plane, we get a plane figure delimited by a 
closed oriented curve; we assume that to 
the latter there corresponds a bi-vector which magnitude gives the planar 
area of the figure, and that the given surface can be decomposed into pieces having the same property.  

Let us decompose the given surface into pieces and, after carrying these pieces arbitrarily in the space, 
let us project all these pieces on the same plane. The sum of magnitudes of the closed oriented curves of  these 
projections  depends on the decomposition of the surface and on the way its pieces are located. 
The supremum of these sums  will be defined as the  {\it bi-vector area} of the surface.
\end{quote}


In 1890 \textsc{Peano} in \emph{Sulla definizione dell'area di una superficie} \cite{peano_area1890} 
examines historically various definitions of area and restates his definition. 
He starts by presenting the definitions of length of a convex planar arc
and the area of a convex surface, given by \textsc{Archimedes}, as the limit 
of inscribed and circumscribed polygons and, respectively, as the limit of inscribed and circumscribed 
convex polyhedral surfaces. 
\textsc{Peano}, aware of the fact that \textsc{Archimede}'s proposal is suitable enough to define the 
area of a cylindrical surface, tries to propose a definition of area preserving the analogy between 
length of arcs and area of surfaces present in \textsc{Archimede}'s work. 
In the case of non planar curves a good definition of length can be obtained by  considering only the inscribed polygons, but in the case of  surfaces, \textsc{Peano} observes that \textsc{Archimede}'s definition cannot be applied to non convex ones.
\textsc{Peano}'s aim is to extend Archimede's definition in order to handle more general 
surfaces, such as the concave ones. 

Later, \textsc{Peano} criticizes the definitions of area present in the literature, including Serret's definition, and explaining that
\begin{quote}
L'errore principale commesso da \textsc{Serret} sta nel ritenere che il piano passante per tre punti di una superficie abbia per limite il piano tangente alla medesima.
\footnote{ \ \pq{The main mistake of \textsc{Serret} is his belief that the plane passing through three points of a surface 
tends to the tangent plane.}}
 \end{quote}
He criticizes also \textsc{Lagrange}'s definition:
\begin{quote}
Il risultato è ottenuto da \textsc{Lagrange} per mezzo di un'asserzione non esatta.
\footnote{ \ \pq{The result has been obtained by \textsc{Lagrange} by means of a not exact statement.}}
\end{quote}
He also criticizes \textsc{Harnach}'s modification of \textsc{Serret}'s definition, saying that, even if the faces of the 
polyhedron considered by \textsc{Harnach} tend to the tangent planes, \textsc{Harnach}'s definition fails even in the case of 
a cartesian surface of equation $z=f(x,y)$.
\textsc{Peano} also recalls that the non correctness of \textsc{Serret}'s definition has already been  noted by \textsc{Schwarz}. 
The definition proposed by \textsc{Hermite} as a consequence of \textsc{Schwarz}'s remark, even if considered sufficiently ``rigorous'' by \textsc{Peano}, is not completely satisfactory, because depends on the choice of the coordinate system.

Finally \textsc{Peano} observes that any difficulty can be overcome by using the concept of {\it oriented} area, 
attributed by him to \textsc{Chelini}, \textsc{M\"obius}, \textsc{Bellavitis}, \textsc{Grassmann} and \textsc{Hamilton}
\footnote{ \ It is interesting to note as \textsc{Chelini}, \textsc{M\"obius}, \textsc{Bellavitis} and \textsc{Grassmann} 
in their work refer directly to \textsc{Poinsot}.}.  
The bi-vectorial definition of non-planar surface area of \textsc{Peano} 
is based on the concept of \textsc{Grassmann}'s {\it bi-vector}: 
\textsc{Peano} extends the equipollence between closed polygonal lines and triangles (see Theorem \ref{2.2}).
Thus closed lines are represented by bi-vectors:
\begin{quotation}
Data una linea chiusa (non piana) $l$, si pu\`o sempre determinare una linea piana chiusa o bivettore $l'$, in guisa che, proiettando le due linee $l$ $l'$ su d'un piano arbitrario, con raggi paralleli di direzione arbitraria, le aree [con segno] limitate dalle loro proiezioni risultino sempre uguali.\footnote{ \ \pq{Given a closed (not planar) line $l$, it is always possible to determine a closed planar line or bi-vector $l'$ in such a way that by projecting both lines on an arbitrary plane, with parallel  rays along an arbitrary direction, the [signed] areas defined by  their projections coincide.}}
\end{quotation}
The logical evidence of this proposition is not trivial for a modern reader 
\footnote{ \ In a letter written by \textsc{Peano} to \textsc{Casorati} (26 October 1889), 
\textsc{Peano} presents his note, later published in 1990 in Rendiconti dell'Accademia dei Lincei:
\begin{quotation}
Questi contorni chiusi sono analoghi ai {\it segmenti} o {\it vettori}, cui corrispondono per dualità; si possono identificare con le {\it coppie} della meccanica. Siccome se\-condo Grassmann sono i prodotti di due vettori, si possono chiamare  {\it bivettori}.\\
Dicasi grandezza di un bivettore $C$ l'area, in valore assoluto, del triangolo $T$ di cui si parla nel precedente teorema. Se proiettando $C$ su una terna di piani ortogonali si ottengono le aree $a,b,c$ allora la grandezza di $C$ vale $\sqrt{a^2+b^2+c^2}$.\\
I bivettori si possono sommare, o comporre, analogamente ai vettori, e precisamente come le coppie di forze. Se una porzione di superficie si scompone in parti, il bivettore (o contorno) di quella superficie è la somma dei bivettori delle sue parti come se un arco di linea si scompone in parti, il vettore (corda) dell'arco è uguale alla somma (risultante) dei vettori delle sue parti.
\end{quotation}
\translation{These closed contours are analogous, by duality, to {\it segments} or {\it vectors}; they can be identified with the couples of mechanics. According to Grassmann they are products of two vectors and can be called {\it bi-vectors}.
Let us call the magnitude of a bi-vector $C$ the area, in absolute value, of the triangle $T$ described in the previous Theorem. 
If, by projecting $C$ on a tern of orthogonal planes, one obtains the areas $a,b,c$, then the size of $C$ is given by $\sqrt{a^2+b^2+c^2}$.
The bi-vectors can be added, or composed, in an analogous way as the vectors, and more precisely as the couples of forces. 
If a part of a surface is decomposed into pieces, the bi-vector (or contour) of that surface is the sum of the bi-vectors of its pieces, 
as in the case of an arc of a line is decomposed into pieces, the vector (cord) of the arc is the sum (resultant) of the vectors of its pieces.}
}.

By presenting the mathematical instruments for the proof, we observe what \textsc{Peano} says in order to understand the necessary mathematical background. 
\begin{quote}
questa proposizione è conseguenza immediata della somma, o composizione, dei bivettori 
[poich\'e tale somma \`e essa stessa un bivettore,] quando la linea $l$ è poligonale.\footnote{ \ \pq{This proposition is a direct consequence of the sum, or composition, of bi-vectors [since such a sum is a bi-vector,] when the line $l$ is polygonal.}}
\end{quote}
The trivialness of this part is a consequence of Theorem \ref{2.2} of Section \ref{sez-vec-calculus}.

\begin{quote}
Il solito passaggio al limite permette di dimostrarla quando la $l$ è una linea curva, descritta da un punto avente sempre derivata finita, ed anche in altri casi.\footnote{ \ \pq{The usual limiting procedure allows one to prove this fact when $l$ is described by a point having finite derivative, and also in other cases.}}
\end{quote}
Concerning this part, the approximation of a line by means of polygons provides the direct way to transfer properties of closed polygons to closed continuous curves. It is worthwhile to note that limits of polygons and triangles are included in the topological concepts introduced by \textsc{Peano} concerning geometric forms 
(see Section \ref{sez-vec-calculus}).     
The condition of finite derivative, besides guaranteeing the continuity of the curve,  assures that any projection of the closed line is the boundary of a set which is measurable in the sense of Jordan-Peano.

Moreover, \textsc{Peano} underlines that area must be thought as ``oriented'':
\begin{quote}
Le aree si devono considerare tenendo in debito conto i segni.
\footnote{ \ \pq{Areas must be considered by taking their sign into account.}}
\end{quote}
This part underlines the fact that the orientation of  closed lines has always to be taken into account and this element becomes fundamental in the case of self-intersecting lines.

Thanks to the notion of equipollence between closed lines, \textsc{Peano} observes that:
\begin{quote}
se si proietta ortogonalmente una linea chiusa (non [necessariamente] piana) $l$ su un piano variabile, il massimo dell'area limitata dalla proiezione di $l$ vale la grandezza del bivettore [associato a] $l$; e questo massimo avviene quando il piano su cui  si proietta ha la giacitura [del bivettore associato a] $l$.
\footnote{ \ \pq{If one projects orthogonally a (not planar) closed line $l$ on a variable plane, the maximum of the area delimited by the projection of $l$ is equal to the size of the bi-vector $l$. This maximum is achieved by projecting on a plane on which $l$ lies.} }
\end{quote}

In 1890 \textsc{Peano} presents a new and more clear formulation of its definition of area:
\begin{quote}
L'area di una porzione di superficie è il limite superiore  della somma delle grandezze dei bivettori delle 
sue parti.\footnote{ \ \pq{The area of a portion of surface is the upper limit of the sum of its parts.}}
\end{quote}


More pragmatically, this quotation suggests the following re-statement of ``bi-vectorial definition of area'':
\begin{quote}
Given an arbitrary non planar surface, we consider a decomposition into pieces. 
For each of these pieces we consider its oriented boundary and the magnitude of the corresponding
bi-vector. 
The supremum, with respect to all decompositions of the surface, of the sums of the magnitudes
of the bi-vectors of the pieces of the decomposition,  
will be defined as the bi-vector area of the surface.
\end{quote}


With this formulation \textsc{Peano} provides the fundamental property leading to the formula of area (\ref{Lagrange}): 
\begin{quote}
La giacitura del  bivettore di una porzione infinitesima di superficie è quella del piano tangente; il rapporto fra la sua grandezza e l'area di quella porzione è l'unità.\footnote{ \ \pq{The bi-vector corresponding to an infinitesimal part of the surface lies on the tangent plane; the ratio between its size and the area of that part is equal to 1.}}
\end{quote}

In this way \textsc{Peano} shows the complete analogy between length of arcs and area of surfaces: in fact \textsc{Peano} 
observes that  the direction of the vector with endpoints on an infinitesimal arc coincides with the tangent, and the rate between their lengths is equal to one. 
\footnote{ \ Besides these properties of areas, \textsc{Peano} gives an estimate of the difference between the lengths 
of an arc and its cord and between the area of a surface and its bi-vector.}
Commenting on this fact, we may say that \textsc{Peano}'s definition grasps the essence of
the measure of area at the infinitesimal level.


The idea  of projection on planes and the selection of the projection which maximizes the area is present also in 
\textsc{Carath\`eodory}'s work of 1914. His ideas are  further developed by \textsc{Hausdorff}, who extends  
\textsc{Carath\`eodory}'s results in the case of \textsc{Hausdorff} measures with integer exponent. Nowadays the most famous measure is the \textsc{Hausdorff} measure, which allows one to define the measure of rather general sets by including 
also the concept of dimension. One of the first results proved by \textsc{Hausdorff} is the \textsc{Lagrange} formula of area \eqref{Lagrange}.

\section{Oriented closed curves and bi-vectors} \label{bi-vectors}

For convenience of the reader we outline in a formal way how to associate a bi-vector to a closed oriented
curve accordingly with Theorem \ref{2.2}.

In \emph{Calcolo geometrico} \cite[(1888)]{peano88} \textsc{Peano} provides a formula to valuate the bi-vector
associated to a closed curve. 
Let $A \colon [t_0, t_1] \to \Pu_3$ be a $C^1$ closed curve. 
The bi-vector associated to $A$ is given by
\footnote{ Recall that $\Pu_3$ denotes the set of points according to Grassmann-Peano vector calculus
(see Section \ref{sez-vec-calculus}) and that $A'(t)$ denotes the derivative of $A$ at $t$.
}

\begin{equation}  \label{int_area}
\int_{t_0}^{t_1} A(t) A'(t) dt  \,\,\,.
\end{equation}
In the case of a closed planar curve $A$, the area of the triangle $X \int_{t_0}^{t_1} A(t) A'(t) dt$
does not depend on the point $X \in \Pu_3$ belonging to the plane of the curve; such area is called by
\textsc{Peano} ``area delimited by the closed planar curve $A$''. 

\textsc{Peano} gives an example of bi-vector associated to the non-planar closed curve
$A: [-r, 2\pi + r + h 2\pi] \to \Pu_3$
formed by a cylindrical helix of radius $r$ and pitch $h 2 \pi$ and three rectilinear pieces,
two horizontal and one vertical, according to the definition:
\begin{equation}
A(t) := \left\{
\begin{array}{l l}
O + (r + t){\bf i} & \text{for  }\   -r \le t \le 0  \\
O+r \cos t{\bf i} + r \sin t {\bf j} + ht {\bf k} & \text{for  }\  0 \le t \le 2\pi  \\
O+(2 \pi + r - t){\bf i} + h 2\pi{\bf k} & \text{for  }\  2\pi \le t \le 2\pi + r  \\
O+(h 2\pi + 2\pi + r -t ){\bf k} & \text{for  }\  2\pi + r \le t \le 2\pi + r + h 2\pi
\end{array}
\right.
\end{equation}
where $O$ is a point of $\Pu_3$ and $\{{\bf i}, {\bf j}, {\bf k}\}$ is an orthogonal base of $\V$.
A straightforward calculations gives for the integral (\ref{int_area}) the value $2 \pi r^2 {\bf ij}$, 
where ${\bf ij}$ denotes the bi-vector product of $\bf i$ and $\bf j$.
This coincides with the value of the bi-vector corresponding to the orthogonal projection of the curve $A$
on the plane ${\bf i},{\bf j}$. 

In previous Section we have outlined several properties related to bi-vectors associated to closed curves.
Now we formulate these properties in terms of the following propositions, leaving the proofs to the reader.  

\begin{proposition}
Let $\gamma \colon [0,1] \to \Pu_3$ be a continuous closed curve lying on a plane $\pi$.
There exists a unique bi-vector (denoted by $\alpha_\gamma$) associated to $\gamma$ such that for any point $O$  
in $\pi$
\begin{equation}
O\alpha_\gamma = \lim_{\{t_i\}} \sum_{i=0}^{m-1}  O \gamma(t_i) \gamma(t_{i+1})
\end{equation}
where the limit is evaluated on the subdivisions $0=t_0 < \dots < t_i < t_{i+1} < \dots < t_m = 1$ of the interval $[0,1]$ 
for $\max \{ t_{i+1} - t_i : i= 0, \dots, m-1 \} \to 0$.
\end{proposition} 

\begin{theorem}
Let $\mu \colon [0,1] \to \Pu_3$ be a continuous closed (non necessarily planar) curve.
There exists a unique bi-vector associated to $\mu$ such that, for every plane $\pi$
and for every parallel projection on $\pi$,
the bi-vector $\alpha_{\mu^*}$ associated to the curve $\mu^*$, projection of the curve $\mu$ 
on the plane $\pi$, is equal to
the projection of the bi-vector on $\pi$.
\end{theorem}

\begin{theorem}
Let $\sigma \colon [0,1]\times [0,1] \to \Pu_3$ be a $C^1$ surface.
For any $x,y \in (0,1)$ let's consider the infinitesimal square 
$Q_\epsilon = [x, x+\varepsilon]\times [y, y+\varepsilon]$,
its counterclockwise oriented boundary $\partial^+ Q_\varepsilon$,
and the infinitesimal element of surface $\sigma(Q_\varepsilon)$.
Then the ratio between the magnitude of the bi-vector $\alpha_{\sigma(\partial^+ Q_\varepsilon)}$,
associated to the closed curve $\partial^+ Q_\varepsilon$,
and the \textsc{Peano}'s area of $\sigma(Q_\varepsilon)$ tends to $1$ when $\varepsilon$ tends to $0$.    
\end{theorem}

\section{Use of the concept of area by \textsc{Peano}} \label{uso}

%
%
%
%
%
In this section we analyze the use of the concept of area by\textsc{Peano} into the following works:
\emph{Applicazioni geometriche} \cite[(1887)]{peano87},
\emph{Calcolo geometrico} \cite[(1888)]{peano88},
\emph{Lezioni di analisi infinitesimale} \cite[(1893)]{peano1893}
and \emph{Formulario mathematico} (1895-1908).

\textsc{Peano}, by means of the notions of inner and outer measures on Euclidean spaces of 
dimension $1, 2, 3\/$, that have been introduced by him in \emph{Sull'integrabilit\`a delle funzioni} \cite[(1883)]{peano1883},
refounds in \emph{Applicazioni geometriche} \cite[(1887)]{peano87} the notion of Riemann integral
and extends it to abstract measures. The development of the theory of measure is based 
on a solid topological and logical ground and on a deep knowledge of set theory.

\textsc{Peano} in \emph{Applicazioni geometriche} and later \textsc{Jordan} in 
\emph{Cours d'Analyse} \cite[(1893)]{jordan1893} develop the well known concepts of classical
measure theory: measurability, change of variables, fundamental theorems of calculus.

The mathematical tools employed by \textsc{Peano} were really innovative both on geometrical and 
topological level. \textsc{Peano} used extensively the geometric vector calculus introduced by
\textsc{Grassmann} (see Section \ref{sez-vec-calculus}). 
A revolutionary tool is the notion of differentiation of distributive set functions,
that suggests to regard area of a non-planar surface as a distributive set function and to compare
it, at the infinitesimal level, with the area of a planar set. In this context the evaluation
of the area of a non-planar surface is reduced to the integration of a numerical function
obtained by differentiation of the area of a non-planar surface with respect to the
area of planar sets.
\footnote{
As observed in our paper \emph{Peano on derivative of measures: strict derivative of distributive
set functions} \cite[(2010)]{der_measures}, differentiation of distributive set functions
gives a mathematical implementation of the \emph{massa-density paradigm} (mass and volume
are distributive set functions and the density is obtained by differentiating mass with
respect to volume). 
}

In this rich mathematical context \textsc{Peano} gives his first definition of area of non-planar surfaces
(see first quotation of Section \ref{sez-Peano-area}) and derives general formulae for planar
and non-planar surfaces.

\begin{enumerate}
\item \labelpag{for-planar}
\emph{Formula for planar area} (see \cite[(1887, Th. 47, p.\,237)]{peano87}, 
\cite[(1893, Vol. 2, \S 394 p.\,224-225)]{peano1893}).
Let $A,B : [t_0,t_1] \to \R^2$ be two $C^1$ functions, such that the segments $A(t)B(t)$ and
$A(t')B(t')$ have empty intersection for any $t, t'  \in [t_0,t_1], t\ne t'$. The set spanned
by the segment $A(t)B(t)$, with $t \in [t_0,t_1]$, namely the set $\cup_{t \in [t_0,t_1]} A(t)B(t)$
has an area $u$ given by the formula
\begin{equation*}  
u = \frac{1}{2} \int_{t_0}^{t_1} (B(t)-A(t)) \cdot \left( \frac{dA(t)}{dt} +  \frac{dB(t)}{dt} \right) dt  
\end{equation*}
where, following \textsc{Peano}'s terminology,
 $(B(t)-A(t)) \cdot \left( \frac{dA(t)}{dt} +  \frac{dB(t)}{dt} \right)$ denotes the
magnitude of the bi-vector given by the product of the vectors $(B(t)-A(t))$ and $ \left( \frac{dA(t)}{dt} +  \frac{dB(t)}{dt} \right)$.
\footnote{ In modern language, the magnitude of this bi-vector is the norm of the vector product
$(B(t)-A(t)) \land \left( \frac{dA(t)}{dt} +  \frac{dB(t)}{dt} \right)$. Therefore
the formula (\ref{for-planar}) becomes
\begin{equation*}
u = \frac{1}{2} \int_{t_0}^{t_1} \left\| (B(t)-A(t)) \land \left( \frac{dA(t)}{dt} +  \frac{dB(t)}{dt} \right) \right\|dt  \,\,\,.
\end{equation*}
}
\item \labelpag{for-non-planar}
\emph{Formula for non-planar area} (see \cite[(1887, Th. 49, p.\,243)]{peano87},
\cite[(1893, Vol. 2, \S 396 p.\,229-232)]{peano1893}).
Let $P: D \to \R^3$ be a $C^1$ function over $D : = \{(u,v) \in \R^2: a<u<b, \,\, \theta_0(u) <v < \theta_1(u)\}$
where $\theta_0$ and $\theta_1$ are continuous functions defined on the interval $[a,b]$.
The surface formed by points $P(u,v)$, with $(u,v)\in D$, has an area $S$ given by the formula
\begin{equation}  \label{bi-area}
S= \int_a^b  du  \int_{\theta_0(u)}^{\theta_1(u)}  \omega (u,v) \, dv
\end{equation}
where $\omega (u,v)$ is the magnitude of the bi-vector product of the vectors $\frac{\partial P}{\partial u}$ and
$\frac{\partial P}{\partial v}$.~%
\footnote{ In modern language, the magnitude of this bi-vector is the norm of the vector product
$\frac{\partial P}{\partial u} (u,v) \land \frac{\partial P}{\partial v} (u,v)$. Therefore
the formula (\ref{for-non-planar}) becomes
\begin{equation*}
S =   \int_a^b  du  \int_{\theta_0(u)}^{\theta_1(u)} \left\|   
\frac{\partial P}{\partial u} (u,v) \land \frac{\partial P}{\partial v} (u,v)   \right\|   \, dv   \,\,\,.
\end{equation*}
}
\end{enumerate} 

\textsc{Peano} uses formulae (\ref{for-planar}) and (\ref{for-non-planar}) to obtain 
classical formulae for elementary surfaces (planar and non-planar). 
Moreover from (\ref{for-planar}) he derives in (\cite[(1887), p.\,242]{peano87}) and in
(\cite[(1893, Vol. 2, \S 394 p.\,225-226)]{peano1893})
formulae that have been recovered
one century later by \textsc{Mamikon A. Mnatskanyan} in his paper
\emph{On the area of the region on a developable surface} \cite[(1981)]{mamikon}.

Particular instances of formula (\ref{for-planar}), considered by \textsc{Peano}, are the following:
\begin{enumerate}
\item The point $A$ moves along a straight line and the angle of the segment $AB$ with that line is constant;

\item The point $A$ is fixed;

\item The segment $AB$ is tangent at the point $A$ to the curve described by $A;$ \labelpag{mamikon-case}

\item The segment $AB$ is of constant length and normal to the curve described by its
midpoint.
\end{enumerate}

In the case (\ref{mamikon-case}), formula (\ref{for-planar}) becomes
\begin{equation*}
u = \frac{1}{2}\int_{t_0}^{t_1}\left\vert \det \left( 
\begin{array}{cc}
v_{1}(t) & v_{2}(t) \\ 
v_{1}^{\prime }(t) & v_{2}^{\prime }(t)%
\end{array}%
\right) \right\vert dt,
\end{equation*}%
where $v_{1}(t),v_{2}(t)$ are the components of $B\left( t\right)
-A\left( t\right) $ and $t\in \left[ t_0,t_1\right] .$ It is clear from this
formula, that the area depends only on the differences of the points $B\left( t\right) - A\left( t\right) $ and not on the particular
positions of the points $A\left( t\right)$, $B\left( t\right)$.
As a consequence of this \textsc{Peano} derives the content of what is nowadays 
stated as \textsc{Mamikon}'s Theorem: 
\emph{the area of a tangent sweep of a curve is equal to the area of its corresponding
tangent cluster}.
The three figures have the same area, because they
are swept by the same tangent vector to the inner ellipsis (or point). The
areas marked by the same letter have the same area as well.
$$
\begin{array}{c}
\includegraphics[scale=0.75]{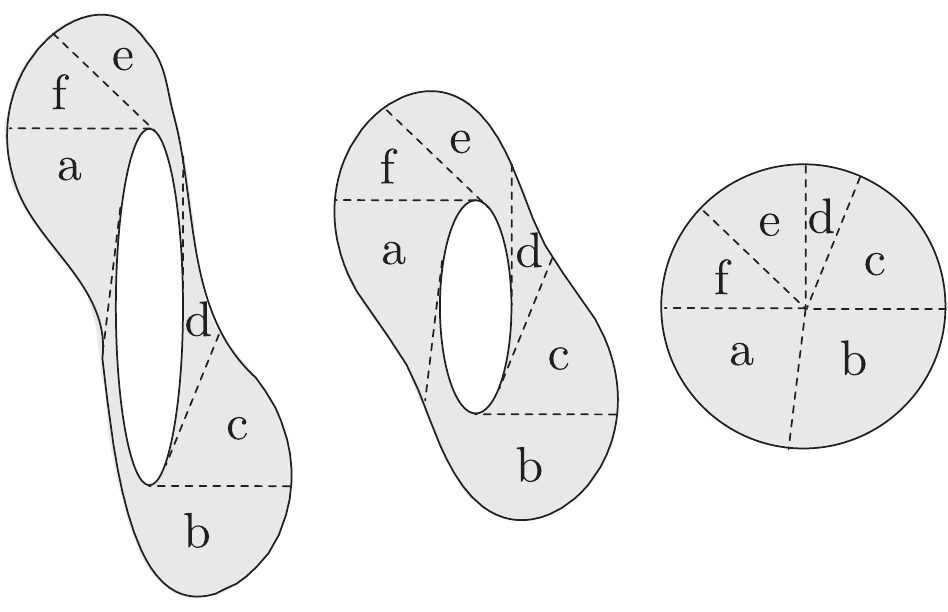}
\end{array}%
$$

\textsc{Mamikon}'s theorem has numerous
applications, as it enables one to obtain area of complicated figures
almost without calculation, by reducing the problem to the calculus of area of
simple figures; see, for examples, \textsc{Mamikon A. Mnatskanyan} and \textsc{Apostol} in
\cite{mamikon-apostol},\cite{mamikon-apostol2}, \cite{mamikon-apostol3}.
\footnote{
In \cite[(2009)]{mamikon-apostol4} \textsc{Apostol} and \textsc{Mnatskanyan}, using
Mamikon theorem, prove the
property of \textsc{Roberval}: ``The area of a cycloidal sector is three times the area
described by the generating disk along its motion''. This property was
proved by \textsc{Peano} \cite[(1893) Vol. 2, \S 395 p.\,226-228)]{peano1893} using (\ref{for-planar}). 
}

Finally the formula (\ref{bi-area}), already obtained by \textsc{Peano} from his geometric definition 
of area of surfaces, is proved by him also using his bi-vectorial definition.
This coincidence is valid in the case of $C^1$ surfaces, but it does not hold for arbitrary
surfaces.

Concerning the area, in the five editions of \emph{Formulario mathematico}, 
in addition to some properties outlined above, we find:
another definition of area \cite[(1902) p. 300]{formulaire4}, due to \textsc{Borchardt}, and
the well-known counter-example to the definition of \textsc{Serret} on area \cite[(1902) p. 300-301]{formulaire4}.

In \emph{Formulario mathematico} \textsc{Peano} adopts
\textsc{Borchardt}'s area \cite[(1854) p. 369]{borchardt}
\footnote{ \textsc{Borchardt}'s area, usually called Minkowski area, was rediscovered by
\textsc{Minkowski} \cite[(1901)]{minkowski} 47 years later.},
defined for every set $S$ of points in $\R^3$ of null volume by the following formula: 

\begin{equation}
\lim_{h \to 0^+} \frac{\text{Volum} \{x \in \R^3 : \text{dist}(x,S) < h  \}}{2 h}
\end{equation}

The counter-example to \textsc{Serret}'s definition is based on the construction
of a polyhedral surface ${\mathcal S}_{m,n}$, with $m,n$ positive integers, 
inscribed into a cylinder of height $1$ and radius $1$,
formed by $mn$ triangles with the following vertices:
\begin{equation*}
\begin{array}{l}
(\cos[\frac{2\pi r}{m}], \sin[\frac{2\pi r}{m}], \frac{s}{n}), 
\, (\cos[\frac{2\pi [r+1]}{m}], \sin[\frac{2\pi [r+1]}{m}], \frac{s}{n}),
\, (\cos[\frac{\pi [2r+1]}{m}], \sin[\frac{\pi [2r+1]}{m}], \frac{s+1}{n})
\end{array} 
\end{equation*}
and by $mn$ triangles with the following vertices:
\begin{equation*}
\begin{array}{l}
(\cos[\frac{2\pi r}{m}], \sin[\frac{2\pi r}{m}], \frac{s}{n}), 
\, (\cos[\frac{\pi [2r-1]}{m}], \sin[\frac{\pi [2r-1]}{m}], \frac{s+1}{n}),
\, (\cos[\frac{\pi [2r+1]}{m}], \sin[\frac{\pi [2r+1]}{m}], \frac{s+1}{n})
\end{array}
\end{equation*}
with $r=0, 1, \dots, m-1$ and $s= 0, 1, \dots, n-1$.

The following pictures show the positions of vertices of triangles in the plane
development of the cylindrical surface (as appears in \textsc{Peano} \cite[(1902) p. 300-301]{formulaire4},
with $m = 5, n=3$)
and the shape of the polyhedral surface ${\mathcal S}_{m,n}$ (as appears in
\textsc{Hermite} \cite[(1883) p. 36]{hermite} with $m=6, n= 10$).

$$
\begin{array}{cc}
\includegraphics[scale=1.2]{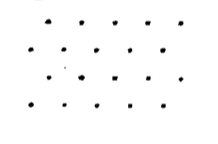} \quad & \quad
\includegraphics[scale=0.1]{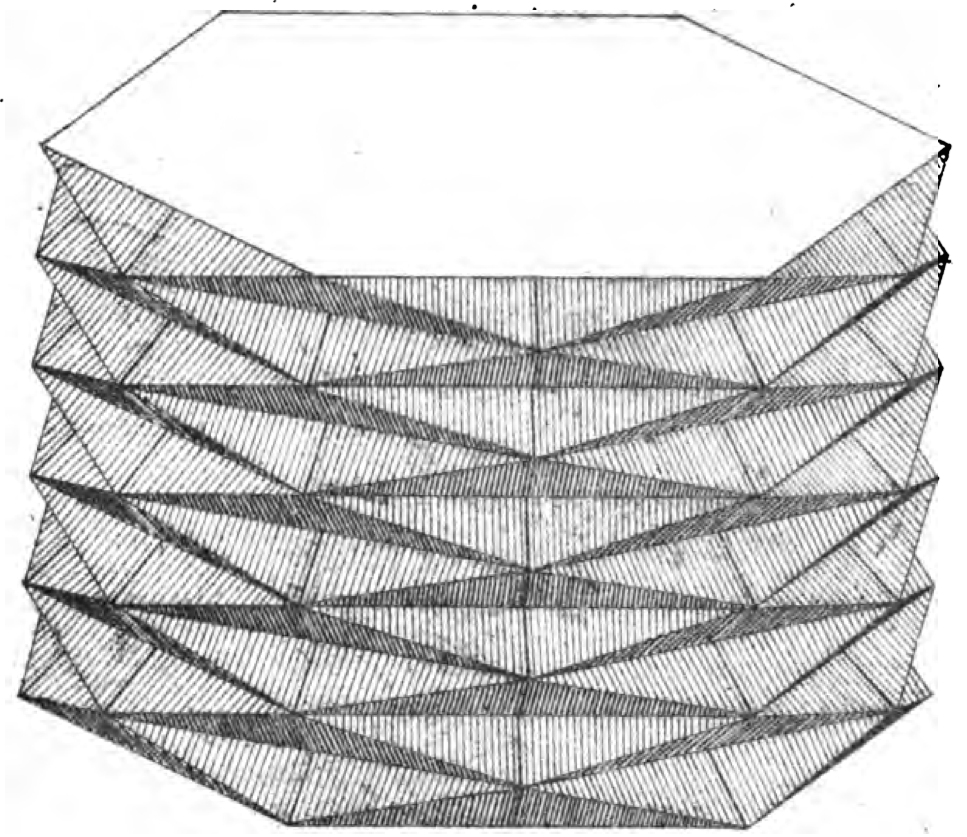}
\end{array}%
$$
\medskip

A straightforward calculations gives the area $a_{m,n}$ of the polyhedral surface ${\mathcal S}_{m,n}$:
\begin{equation}
a_{m,n} = 2m \sin(\frac{\pi}{m}) \sqrt{1 + 4 n^2 \sin^4 \frac{\pi}{2m}}  \quad .
\end{equation}

Clearly
\begin{equation}
\lim_{m \to \infty} a_{m,m} = 2\pi \,, \quad
\lim_{m \to \infty} a_{m,m^2} = 2\pi \sqrt{1 + \frac{\pi^4}{4}} \,, \quad
\lim_{m \to \infty} a_{m,m^3} = + \infty
\end{equation}

Consequently the limit of the area of the polyhedra ${\mathcal S}_{m,n}$ for $m,n \to \infty$ does not exist.

%
%
%
%
%

\section{On the influence of \textsc{Peano} on definition of area} \label{Peano-influenza}

With \textsc{Lebesgue}'s Thesis \emph{Int\'egrale, Longueur, Aire} \cite[(1902)]{lebesgue1902}, 
\textsc{Peano}'s definition of area acquires notoriety. \textsc{Lebesgue}
is acquainted with the bi-vectorial definition of area given by \textsc{Peano} in 1890,
but ignores the original definition of 1887 and any other contribution of this Author (with the
exception of the \textsc{Peano}'s curve). As a consequence of this, it is not surprising that, 
in almost all contributions on the definition of area, 
references to the other \textsc{Peano}'s works on area and, in particular, to the books 
\emph{Applicazioni geometriche} \cite[(1887)]{peano87} 
and \emph{Calcolo geometrico} \cite[(1888)]{peano88}, are absent.

\textsc{Lebesgue}'s area of a parameterized surface is defined by him as the
lower limit of the area of the polyhedral surfaces that approximate uniformly the surface. 

%
%
%
%
%
In the mathematical literature, we find definitions of area that implement 
\textsc{Peano}'s inequality, namely the ``area of surface is greater or equal to the area of
its orthogonal projection on an arbitrary plane''. Different implementations 
correspond to the different way of defining the ``area of the orthogonal projection
on a plane''. 

Other definitions of area implement the \textsc{Peano}'s bi-vectorial inequality, 
namely that the ``area of a surface bounded by a closed oriented contour is greater or equal to the
magnitude of the bi-vector associated with the contour itself''. In this case, the implementations
correspond to the different ways to associate a number to a given oriented closed curve.

After \textsc{Schwarz} and \textsc{Peano}, as observed by \textsc{Rad\'o} in \cite[(1956) p.\,513]{rado-cesari}, 
``many definitions [of surface area] have been proposed, and an enormous amount of efforts
have been expended in the study of \dots various concepts of surface area''.
For this reason we are forced to present only some contributions. Interested readers
may find detailed historical and mathematical facts in \textsc{Cesari}'s \emph{Surface area} \cite[(1954)]{cesari-book}
and \textsc{Rad\'o}'s \emph{Length and area} \cite[(1948]{rado-book}. 

In addition to the one given by \textsc{Peano}, remarkable definitions are the \textsc{Lebesgue}'s and \textsc{Ge\"ocze}'s area.
The original definitions of \textsc{Peano} and \textsc{Ge\"ocze} provide an evaluation of area that is
greater than or equal to \textsc{Lebesgue}'s area. Observe that \textsc{Peano}'s and \textsc{Ge\"ocze}'s area
relies on the evaluation of the area of the orthogonal projection on planes
of pieces of the given surface. Therefore many authors have proposed different  
ways to define the area of a plane surface, in order to make \textsc{Peano}'s and \textsc{Ge\"ocze}'s area
coincident with \textsc{Lebesgue}'s area for a wider class of continuous parametric surfaces
(see \textsc{Rad\'o} \cite[(1928)]{rado1928} and \textsc{Cecconi} \cite[(1950)]{cecconi1950},
\cite[(1951)]{cecconi1951}).  

\textsc{Cesari} \cite[(1956)]{cesari-book} reformulates
the definitions given by \textsc{Peano} and \textsc{Ge\"ocze}
 in a suitable way in order ``to preserve'' elementary area of polyhedral surfaces
and, above all, lower semicontinuity.
\textsc{Cesari} states the following theorem:
\begin{theorem} \label{LVP}
For every continuous surface $S$ we have $\Lc(S) = \Vc(S) = \Pc(S)$,
where $\Lc(S)$, $\Vc(S)$ and $\Pc(S)$ denote Lebesgue area, Ge\"ocze area and Peano area.
\end{theorem}

More precisely, \textsc{Peano}'s and \textsc{Ge\"ocze}'s definitions  are reformulated by \textsc{Cesari} in terms of 
\emph{topological index} of planar closed curves. This index is denoted with $O(P, \gamma)$ by \textsc{Cesari},
where $\gamma$ is a closed planar curve, and $P=(x,y)$ is a point of the plane $\pi$ of $\gamma$. 
It is worth observing that, as for the bi-vector associated to a closed planar curve,
the integral $\int_{\pi} |O(P, \gamma)|  dx dy $ (denoted in the following by $v(\gamma, \pi)$) is interpreted, as ``area of the planar
surface delimited by $\gamma$ '' (see \textsc{Cesari} \cite[(1956) p.\,104]{rado-cesari}).

Now, let $S$ be a parametric surface in $\R^3$, parameterized by a continuous $\varphi\colon A \to S$ (i.e. $S = \varphi (A)$),
where $A$ is an \emph{admissible} set
(\footnote{ Among the admissible sets (see \textsc{Cesari} \cite[(1956) p. 27]{cesari-book}), we mention:
planar sets delimited by a Jordan simple curve or finite union of such sets, and open sets.}).  
Given a plane $\alpha$ and a curve $\gamma$ in $A$, 
let denote with $\gamma^{*\alpha}$ the orthogonal projection
on $\alpha$ of the image $\gamma^*$ of $\gamma$ under the parameterization $\varphi$. 

The reformulation $\Pc(S)$ of \textsc{Peano}'s area of the surface $S$, given by \textsc{Cesari} 
(see \cite[(1956) p.\,137]{cesari-book}), is the following:
\begin{equation} \label{peano-cesari}
\Pc (S) := \sup_{\{ \gamma_i\}_i} \sum_{i} \sup_{\alpha} \, v(\gamma_i^{*\alpha}, \alpha)
\end{equation} 
where $\{ \gamma_i\}_i$ runs over all finite families of simple closed polygonal curves in $A$
delimiting non-overlapping regions and $\alpha$ runs over all planes in $\R^3$.

Concerning \textsc{Ge\"ocze}'s area, let us consider the coordinate planes $\alpha_{xy}$, $\alpha_{yz}$ 
and $\alpha_{zx}$ in the Euclidean space.
The reformulation $\Vc(S)$ of \textsc{Ge\"ocze}'s area of the surface $S$, given by \textsc{Cesari}
(see \cite[(1956) p.\,117]{cesari-book}), is the following:
\begin{equation}  \label{geocze-cesari}
\Vc(S) := \sup_{\{ \gamma_i\}_i} \sum_{i} \sqrt{ [v(\gamma_i^{*\alpha_{xy}}, \alpha_{xy})]^2
+ [(v(\gamma_i^{*\alpha_{yz}}, \alpha_{yz})]^2  +  [(v(\gamma_i^{*\alpha_{zx}}, \alpha_{zx})]^2} 
\end{equation} 

A great deal of research has been dedicated to find an axiomatic characterization
of a notion of surface area, namely to the problem of establishing properties characterizing 
univocally the notion of area. \textsc{Cecconi} in \cite[(1951)]{cecconi1951}
gives the following properties characterizing \textsc{Lebesgue}'s area (and, consequently, Peano area \eqref{peano-cesari}
and Ge\"ocze area \eqref{geocze-cesari}): 
\begin{theorem}
Let $\Phi$ be a functional defined over all continuous parametric surfaces $S$ on $2$-cells.
Then $\Phi$ coincides with Lebesgue area if the following properties are satisfied:
\begin{enumerate}
\item \labelpag{cecc1}
$\Phi$ is lower semi-continuous;
\item \labelpag{cecc2}
$\Phi$ coincides with usual elementary area for polyhedral surfaces;
\item \labelpag{cecc3}
$\Phi$ is super-additive
(\footnote{ Namely, for every subdivision of the surface $S$ in parts,
the area of $S$ is greater than the sum of the areas of the parts.});
\item \labelpag{cecc4}
$\Phi$ satisfies Peano inequality
(\footnote{ Namely, for every $S$, and for every plane $\alpha$, one has
$\Phi(S) \ge \textrm{mis} (\{ P \in \alpha : O(P;C_\alpha) \ne 0  \})$,
where $C_\alpha$ is the projection of the contour of $S$ on $\alpha$, and
$O(P;C_\alpha)$ is the topological index of $P$ with respect to the curve $C_\alpha$ defined above.}).
\end{enumerate}
\end{theorem}
In the proof of this Theorem, given by \textsc{Cecconi}, a crucial step consists in the inequality 
$\Pc(S) \le \Phi(S) \le \Lc(S)$ that, together with the equality $\Pc(S) = \Lc(S)$ (see Theorem \ref{LVP}),
leads to the expected coincidence $\Phi(S) = \Lc(S)$.

%
%
%
%
%
%
\bibliographystyle{plain}

\end{document}